\documentclass[a4paper,11pt]{amsart}

\usepackage{graphicx}
\usepackage{mathptmx}
\usepackage{amsmath}
\usepackage{amssymb}
\usepackage{enumitem}
\usepackage{xcolor}

\newmuskip\pFqmuskip

\newcommand*\pFq[6][8]{%
  \begingroup 
  \pFqmuskip=#1mu\relax
  \mathcode`=\string"8000
  \begingroup\lccode`\~=`\,
  \lowercase{\endgroup\let~}\pFqcomma
  F^{#2}_{#3}{\left(\genfrac..{0pt}{}{#4}{#5}\bigg|#6\right)}%
  \endgroup
}
\newcommand{\pFqcomma}{\mskip\pFqmuskip}

\newtheorem{theorem}{Theorem}[section]

\begin{document}

\title[]{Some identities on degenerate trigonometric functions}

\author{Taekyun  Kim}
\address{Department of Mathematics, Kwangwoon University, Seoul 139-701, Republic of Korea}
\email{tkkim@kw.ac.kr}
\author{Dae San  Kim}
\address{Department of Mathematics, Sogang University, Seoul 121-742, Republic of Korea}
\email{dskim@sogang.ac.kr}


\begin{abstract}
In this paper, we study several degenerate trigonometric functions, which are degenerate versions of the ordinary trigonometric functions, and derive some identities among such functions by using elementary methods. Especially, we obtain multiple angle formulas for the degenerate cotangent and degenerate sine functions.
\end{abstract}

\maketitle

\markboth{\centerline{\scriptsize Some identities on degenerate trigonometric functions}}
{\centerline{\scriptsize T. Kim and D. S. Kim}}

\section{Introduction}
For any nonzero $\lambda\in\mathbb{R}$, the degenerate exponentials are defined by
\begin{equation*}
e_{\lambda}^{x}(t)=\sum_{n=0}^{\infty}(x)_{n,\lambda}\frac {t^{n}}{n!},\quad (n\ge 0),\quad (\mathrm{see}\ [4,5]),
\end{equation*}
where
\begin{displaymath}
(x)_{0,\lambda}=1,\ (x)_{n,\lambda}=x(x-\lambda)(x-2\lambda)\cdots\big(x-(n-1)\lambda\big),\ (n\ge 1).
\end{displaymath}
Note that $\displaystyle\lim_{\lambda\rightarrow 0}e_{\lambda}^{x}(t)=e^{xt}\displaystyle$ and $e_{\lambda}(t)=e_{\lambda}^{1}(t)$. In [4], the degenerate hyperbolic functions are introduced by
\begin{equation*}
\begin{aligned}
&\cosh_{\lambda}\big(x:a\big)=\frac{1}{2}\big(e_{\lambda}^{x}(a)+e_{\lambda}^{-x}(a)\big),\ \sinh_{\lambda}\big(x:a\big)=\frac{1}{2}\big(e_{\lambda}^{x}(a)-e_{\lambda}^{-x}(a)\big), \\
&\tanh_{\lambda}\big(x:a\big)=\frac{\sinh_{\lambda}(x:a)}{\cosh_{\lambda}(x:a)},\quad \coth_{\lambda}(x:a)=\frac{\cosh_{\lambda}(x:a)}{\sinh_{\lambda}(x:a)},\ (x\ne 0,\ a\ne 0).
\end{aligned}	
\end{equation*}
It is well known that
\begin{equation}
\cos a=\frac{e^{ia}+e^{-ia}}{2},\quad \sin a=\frac{e^{ia}-e^{-ia}}{2i}, \label{1}	
\end{equation}
where $i=\sqrt{-1}$, (see [1,2,3]). \par
In light of \eqref{1}, we consider the degenerate cosine and degenerate sine functions which are respectively given by
\begin{equation}
\cos_{\lambda}\big(x:a\big)=\frac{e_{\lambda}^{xi}(a)+e_{\lambda}^{-xi}(a)}{2},\quad \sin_{\lambda}\big(x:a\big)=\frac{e_{\lambda}^{xi}(a)-e_{\lambda}^{-xi}(a)}{2i}. \label{2}	
\end{equation}
Note that
\begin{displaymath}
	\lim_{\lambda\rightarrow 0}\cos_{\lambda}\big(x:a\big)=\cos ax,\quad \lim_{\lambda\rightarrow 0}\sin_{\lambda}\big(x:a\big)=\sin ax.
\end{displaymath}
In addition, we define the degenerate tangent and degenerate cotangent functions as follows:
\begin{equation}
\tan_{\lambda}\big(x:a\big)=\frac{\sin_{\lambda}(x:a)}{\cos_{\lambda}(x:a)},\quad \cot_{\lambda}\big(x:a\big)=\frac{\cos_{\lambda}(x:a)}{\sin_{\lambda}(x:a)},\label{3}	
\end{equation}
where $x\ne 0$ and $a\ne 0$. \par
Recently, degenerate versions of many special numbers and polynomials and some transcendental functions have been investigated. The aim of this paper is to study several degenerate trigonometric functions, namely the degenerate sine, degenerate cosine, degenerate tangent and degenerate cotangent functions. We derive some identities on those degenerate trigonometric functions by using elementary methods. In particular, we show the following multiple angle fomulas for the degenerate cotangent and degenerate sine functions (see \eqref{11}, \eqref{15}, \eqref{28}):

\begin{align*}
&2m\cot_{\lambda}\big(2mt:a\big)=-\sum_{j=0}^{2m-1}\tan_{\lambda}\bigg(t+\frac{j\pi}{2m \log e_{\lambda}(a)}:a\bigg), \\
&2m\cot_{\lambda}\big(2mt:a\big)=\sum_{j=0}^{2m-1}\cot_{\lambda}\Big(t+\frac{j\pi}{2m \log e_{\lambda}(a)}:a\Big), \\
&\sin_{\lambda}\big((2m+1)x:a\big)=(2m+1)\sin_{\lambda}\big(x:a\big) \prod_{k=1}^{m}\Bigg(1-\cfrac{\sin_{\lambda}^{2}(x:a)}{\sin^{2}\Big(\frac{k \pi}{2m+1} \Big)}\Bigg).
\end{align*}

\section{Some identities on degenerate trigonometric functions}
From \eqref{2}, we note that
\begin{equation}
\begin{aligned}
\sin_{\lambda}^{2}\big(x:a\big)+\cos_{\lambda}^{2}\big(x:a\big)&=\bigg(\frac{e_{\lambda}^{xi}(a)+e_{\lambda}^{-xi}(a)}{2}\bigg)^{2}+ \bigg(\frac{e_{\lambda}^{xi}(a)-e_{\lambda}^{-xi}(a)}{2i}\bigg)^{2}\\
&=\frac{4}{4}=1.
\end{aligned}\label{4}
\end{equation}
We observe that
\begin{equation}
\begin{aligned}
\cos_{\lambda}\big(2x:a\big)&=\frac{e_{\lambda}^{2xi}(a)+e_{\lambda}^{-2xi}(a)}{2}=1+\frac{e_{\lambda}^{2xi}(a)+e_{\lambda}^{-2xi}(a)-2}{2} \\
&=1-2\bigg(\frac{e_{\lambda}^{xi}(a)-e_{\lambda}^{-xi}(a)}{2i}\bigg)^{2}=1-2\sin_{\lambda}^{2}\big(x:a\big).
\end{aligned}\label{5}
\end{equation}
By \eqref{4} and \eqref{5}, we get
\begin{equation*}
\cos_{\lambda}\big(2x:a\big)=1-2\sin_{\lambda}^{2}(x:a)=2\cos_{\lambda}^{2}(x:a)-1.
\end{equation*}
Moreover, by \eqref{2}, we get
\begin{equation*}
\begin{aligned}
\sin_{\lambda}\big(2x:a\big)&=\frac{e_{\lambda}^{2xi}(a)-e_{\lambda}^{-2xi}(a)}{2i}=2\frac{e_{\lambda}^{xi}(a)-e_{\lambda}^{-ix}(a)}{2i}\frac{e_{\lambda}^{xi}(a)+e_{\lambda}^{-xi}(a)}{2} \\
&=2\sin_{\lambda}\big(x:a\big)\cos_{\lambda}\big(x:a\big).
\end{aligned}
\end{equation*}
\begin{theorem}
For $a\in\mathbb{R}$, we have
\begin{displaymath}
\cos_{\lambda}\big(2x:a\big)=1-2\sin_{\lambda}^{2}\big(x:a\big)=2\cos_{\lambda}^{2}\big(x:a\big)-1,
\end{displaymath}
and
\begin{displaymath}
\sin_{\lambda}\big(2x:a\big)=2\sin_{\lambda}\big(x:a\big)\cos_{\lambda}\big(x:a\big).
\end{displaymath}
\end{theorem}
\begin{align*}
&\sin_{\lambda}\big(x+y:a\big)=\frac{1}{2i}\Big(e_{\lambda}^{(x+y)i}-e_{\lambda}^{-(x+y)i}(a)\Big) \\
&=\frac{e_{\lambda}^{xi}(a)-e_{\lambda}^{-xi}(a)}{2i}\frac{e_{\lambda}^{yi}(a)+e_{\lambda}^{-yi}(a)}{2}+\frac{e_{\lambda}^{xi}(a)+e_{\lambda}^{-xi}(a)}{2}\frac{e_{\lambda}^{yi}(a)-e_{\lambda}^{-yi}(a)}{2i}\\
&=\sin_{\lambda}\big(x:a\big)\cos_{\lambda}\big(y:a\big)+\cos_{\lambda}\big(x:a\big)\sin_{\lambda}\big(y:a\big),
\end{align*}
and
\begin{align*}
&\cos_{\lambda}\big(x+y:a\big)=\frac{e_{\lambda}^{(x+y)i}(a)+e_{\lambda}^{-(x+y)i}(a)}{2} \\
&=\frac{e_{\lambda}^{xi}(a)+e_{\lambda}^{-xi}(a)}{2}\frac{e_{\lambda}^{yi}(a)+e_{\lambda}^{-yi}(a)}{2}-\frac{e_{\lambda}^{xi}(a)-e_{\lambda}^{-xi}}{2i}\frac{e_{\lambda}^{yi}(a)-e_{\lambda}^{-yi}(a)}{2i} \\
&=\cos_{\lambda}\big(x:a\big)\cos_{\lambda}\big(y:a\big)-\sin_{\lambda}\big(x:a\big)\sin_{\lambda}\big(y:a\big).
\end{align*}
Therefore, we obtain the following theorem.
\begin{theorem}
For $a\in\mathbb{R}$, we have
\begin{align*}
&\sin_{\lambda}\big(x\pm y:a\big)=\sin_{\lambda}\big(x:a\big)\cos_{\lambda}\big(y:a\big)\pm \cos_{\lambda}\big(x:a\big)\sin_{\lambda}\big(y:a\big), \\
&\cos_{\lambda}\big(x\pm y:a\big)=\cos_{\lambda}\big(x:a\big)\cos_{\lambda}\big(y:a\big)\mp\sin_{\lambda}\big(x:a\big)\sin_{\lambda}\big(y:a\big).
\end{align*}
\end{theorem}
We observe that
\begin{align}
\frac{d}{dx}\cos_{\lambda}\big(x:a\big)=-\log e_{\lambda}(a)\sin_{\lambda}(x:a),  \label{6}\\
\frac{d}{dx}\sin_{\lambda}\big(x:a\big)=\log e_{\lambda}(a)\cos_{\lambda}(x:a). \label{7}
\end{align}
From \eqref{2} and noting that $\prod_{j=0}^{2m-1}(x+e^{\frac{2 \pi ij}{2m}})=x^{2m}-1$, we get
\begin{align}
&\prod_{j=0}^{2m-1}\cos_{\lambda}\bigg(t+\frac{j\pi}{2m\log e_{\lambda}(a)}:a\bigg)=\prod_{j=0}^{2m-1}\bigg(\frac{e_{\lambda}^{(t+\frac{j\pi}{2m\log e_{\lambda}(a)})i}(a)+e_{\lambda}^{-(t+\frac{j\pi}{2m \log e_{\lambda (a)}})i}(a)}{2}\bigg) \label{8}\\
&=\frac{1}{4^{m}}\prod_{j=0}^{2m-1}\Big(e_{\lambda}^{it}(a)e_{\lambda}^{-\frac{j\pi}{2m \log e_{\lambda}(a)}i}(a)\Big(e_{\lambda}^{\frac{2j\pi}{2m \log e_{\lambda} (a)}i}(a)+e_{\lambda}^{-2it}(a)\Big)\bigg)\nonumber \\
&=\frac{1}{4^{m}}e^{-\frac{\pi}{2m}i\frac{2m(2m-1)}{2}} e_{\lambda}^{2mti}(a)\prod_{j=0}^{2m-1}\bigg(e_{\lambda}^{-2it}(a)+e^{\frac{2\pi ij}{2m}}\bigg) \nonumber \\
&=\frac{1}{4^{m}}e_{\lambda}^{2mti}(a)e^{-m\pi i} e^{\frac{\pi}{2}i}\Big(e_{\lambda}^{-4mti}(a)-1\Big) \nonumber \\
&=\frac{1}{4^{m}}e^{-m\pi i}e^{\frac{\pi}{2}i}\Big(e_{\lambda}^{-2mti}(a)-e_{\lambda}^{2mti}(a)\Big) \nonumber \\
&=-\frac{1}{4^{m}}e^{-m\pi i}e^{\frac{\pi}{2}i}2i\bigg(\frac{e_{\lambda}^{2mti}(a)-e_{\lambda}^{-2mti}(a)}{2i}\bigg) \nonumber \\
&=-\frac{i}{2^{2m-1}}e^{-m\pi i}e^{\frac{\pi}{2}i}\sin_{\lambda}(2mt:a). \nonumber
\end{align}
Note that
\begin{displaymath}
\bigg|-\frac{i}{2^{2m-1}}e^{-m\pi i}e^{\frac{\pi}{2}i}\sin_{\lambda}\big(2mt:a\big)\bigg|=\bigg|\frac{1}{2^{2m-1}}\sin_{\lambda}\big(2mt:a\big)\bigg|.
\end{displaymath}
Thus, by \eqref{8}, we get
\begin{equation}
\log\bigg|\frac{1}{2^{2m-1}}\sin_{\lambda}\big(2mt:a\big)\bigg|=\sum_{j=0}^{2m-1}\log\bigg|\cos_{\lambda}\bigg(t+\frac{j\pi}{2m \log e_{\lambda}(a)}:a\bigg)\bigg|. \label{9}
\end{equation}
By using \eqref{6}, we take the derivatives with respect to $t$ on both sides of \eqref{9}  and get the following identity:
\begin{equation}
\frac{2m\log(e_{\lambda}(a))\cos_{\lambda}(2mt:a)}{\sin_{\lambda}(2mt:a)}=-\sum_{j=0}^{2m-1}\log\big(e_{\lambda}(a)\big)\frac{\sin_{\lambda}\big(t+\frac{j\pi}{2m \log e_{\lambda}(a)}:a\big)}{\cos_{\lambda}\big(t+\frac{j\pi}{2m \log e_{\lambda}(a)}:a\big)}.\label{10}
\end{equation}
Therefore, by \eqref{10}, we obtain the following theorem.
\begin{theorem}
For $m\in\mathbb{N}$, we have
\begin{equation}
-2m\cot_{\lambda}\big(2mt:a\big)=\sum_{j=0}^{2m-1}\tan_{\lambda}\bigg(t+\frac{j\pi}{2m \log e_{\lambda}(a)}:a\bigg). \label{11}	
\end{equation}
\end{theorem}
We remark here that
\begin{equation*}
\tan_{\lambda} \Big(t+\frac{j\pi}{2m \log e_{\lambda}(a)}:a\Big)=\tan \Big(t \log e_{\lambda} (a)+\frac{j \pi}{2m} \Big).
\end{equation*}
Taking $\lambda \rightarrow 0$ in \eqref{11}, we obtain
\begin{align*}
-2m\cot\big(2mat\big)=\sum_{j=0}^{2m-1}\tan \Big(at+\frac{j \pi}{2m}\Big),\ (\mathrm{see}\ [1,2]).
\end{align*}
Working with $\prod_{j=0}^{2m-1} \sin_{\lambda}\Big(t+\frac{j \pi}{2m \log e_{\lambda} (a)}:a\Big)$ and proceeding analogously to \eqref{8}, we obtain
\begin{equation}
\prod_{j=0}^{2m-1} \sin_{\lambda}\Big(t+\frac{j \pi}{2m \log e_{\lambda} (a)}:a\Big)=-\frac{i}{2^{2m-1}}(-1)^{m}e^{-m\pi i}e^{\frac{\pi}{2}i}\sin_{\lambda}(2mt:a). \label{12}
\end{equation}
Thus, from \eqref{12}, we have
\begin{equation}
\log\bigg|\frac{1}{2^{2m-1}}\sin_{\lambda}\big(2mt:a\big)\bigg|=\sum_{j=0}^{2m-1}\log\bigg|\sin_{\lambda}\Big(t+\frac{j\pi}{2m \log e_{\lambda}(a)}:a\Big)\bigg|. \label{13}
\end{equation}
Using\eqref{7}, we take the derivatives with respect to $t$ on both sides of \eqref{13} and get the identity:
\begin{equation}
\frac{2m\log(e_{\lambda}(a))\cos_{\lambda}(2mt:a)}{\sin_{\lambda}(2mt:a)}=\sum_{j=0}^{2m-1}\log\big(e_{\lambda}(a)\big)\frac{\cos_{\lambda}\big(t+\frac{j\pi}{2m \log e_{\lambda}(a)}:a\big)}{\sin_{\lambda}\big(t+\frac{j\pi}{2m \log e_{\lambda}(a)}:a\big)}.\label{14}
\end{equation}

Therefore, by \eqref{14}, we obtain the following theorem.
\begin{theorem}
For $m\in\mathbb{N}$, we have
\begin{equation}
2m\cot_{\lambda}\big(2mt:a\big)=\sum_{j=0}^{2m-1}\cot_{\lambda}\Big(t+\frac{j\pi}{2m \log e_{\lambda}(a)}:a\Big). \label{15}	
\end{equation}
\end{theorem}
Taking $\lambda \rightarrow 0$ in \eqref{15}, we obtain
\begin{align*}
2m\cot\big(2mat\big)=\sum_{j=0}^{2m-1}\cot(at+\frac{j \pi}{2m}),\ (\mathrm{see}\ [1,2]).
\end{align*}

Now, we observe that
\begin{equation}
\cos_{\lambda}\big((k+1)x:a\big)+\cos_{\lambda}\big((k-1)x:a\big)=2\cos_{\lambda}\big(kx:a\big)\cos_{\lambda}\big(x:a\big).\label{16}	
\end{equation}
Thus, by \eqref{16}, we get
\begin{align}
&\cos_{\lambda}\big((k+1)x:a\big)=2\cos_{\lambda}\big(kx:a\big)\cos_{\lambda}(x:a)-\cos_{\lambda}\big((k-1)x:a\big)\label{17}\\
&=2^{2}\cos^{2}_{\lambda}\big(x:a\big)\cos_{\lambda}\big((k-1)x:a\big)-2\cos_{\lambda}\big(x:a\big)\cos_{\lambda}\big((k-2)x:a\big)-\cos_{\lambda}\big((k-1)x:a\big) \nonumber\\
&\cdots\nonumber\\
&=2^{k+1}\cos_{\lambda}^{k+1}\big(x:a\big)+\cdots \nonumber
\end{align}
From \eqref{17}, we note that there is polynomial $T_{n}(x)$ with degree $n$ such that
\begin{equation}
\cos_{\lambda}\big(nx:a\big)=T_{n}\Big(\cos_{\lambda}\big(x:a\big)\Big),\quad (n\in\mathbb{N}). \label{18}
\end{equation}
By using Theorem 2.2, we get
\begin{equation}
\sin_{\lambda}\big((2k+1)x:a\big)-\sin_{\lambda}\big((2k-1)x:a\big)=2\cos_{\lambda}\big(2kx:a\big)\sin_{\lambda}\big(x:a\big).\label{19}
\end{equation}
By \eqref{18} and Theorem 2.1, we get
\begin{equation}
\cos_{\lambda}\big(2kx:a\big)=T_{k}\Big(\cos_{\lambda}\big(2x:a\big)\Big)=T_{k}\Big(1-2\sin_{\lambda}^{2}\big(x:a\big)\Big). \label{21}
\end{equation}
From \eqref{18}, \eqref{19} and \eqref{21}, we show that
\begin{equation*}
\sin_{\lambda}\big((2m+1)x:a \big)
=\sin_{\lambda}\big(x:a \big)\bigg(1+2\sum_{k=1}^{m}T_{k}\Big(1-2\sin_{\lambda}^{2} \big(x:a \big)\Big)\bigg),
\end{equation*}
and hence there exist polynomials $K_{m}(x)$ with degree $m$ such that
\begin{equation}
\sin_{\lambda}\Big((2m+1)x:a\Big)=\sin_{\lambda}(x:a)K_{m}\Big(\sin_{\lambda}^{2}(x:a)\Big).\label{22}
\end{equation}
Now, from \eqref{22} we note that
\begin{equation}
\begin{aligned}
0&=\sin_{\lambda}\bigg((2m+1)\frac{k\pi}{(2m+1)\log(e_{\lambda}(a))}:a\bigg) \\
&=\sin_{\lambda}\bigg(\frac{k\pi}{(2m+1)\log(e_{\lambda}(a))}:a\bigg)K_{m}\bigg(\sin_{\lambda}^{2}\bigg(\frac{k\pi}{(2m+1)\log e_{\lambda}(a)}:a\bigg)\bigg) \\
&=\sin \Big(\frac{k \pi}{2m+1} \Big)K_{m}\Big(\sin^{2}\Big(\frac{k \pi}{2m+1} \Big)\Big),
\end{aligned}\label{23}
\end{equation}
which implies that $K_{m}\Big(\sin^{2}\Big(\frac{k \pi}{2m+1} \Big)\Big)=0$, for $k=1,2,3,\dots,m$. \par
By \eqref{23}, the polynomials $K_{m}(x)$ can be written as
\begin{equation}
K_{m}(x)=C\prod_{k=1}^{m}\bigg(1-\cfrac{x}{\sin^{2}\big(\frac{k \pi}{2m+1} \big)}\bigg).\label{24}
\end{equation}
Note from \eqref{22} that
\begin{align}
C&=K_{m}(0)=\lim_{x\rightarrow 0}K_{m}\big(\sin_{\lambda}^{2}(x:a)\big)\label{25} \\
&=\lim_{x\rightarrow 0}\frac{\sin_{\lambda}\big((2m+1)x:a\big)}{\sin_{\lambda}(x:a)}=\frac{(2m+1)\log(e_{\lambda}(a))}{\log e_{\lambda}(a)}=2m+1. \nonumber
\end{align}
By \eqref{24} and \eqref{25}, we get
\begin{equation}
K_{m}(x)=(2m+1)\prod_{k=1}^{m}\bigg(1-\cfrac{x}{\sin^{2}\big(\frac{k \pi}{2m+1} \big)}\bigg).\label{26}	
\end{equation}
From \eqref{22} and \eqref{26}, we have
\begin{align}
&\sin_{\lambda}\big((2m+1)x:a\big)=\sin_{\lambda}\big(x:a\big)K_{m}\big(\sin_{\lambda}^{2}(x:a)\big)\label{27} \\
&=(2m+1)\sin_{\lambda}\big(x:a\big) \prod_{k=1}^{m}\bigg(1-\cfrac{\sin_{\lambda}^{2}(x:a)}{\sin^{2}\big(\frac{k \pi}{2m+1} \big)}\bigg).\nonumber
\end{align}
Therefore, by \eqref{27}, we obtain the following theorem.
\begin{theorem}
For $m\in\mathbb{N}$, we have
\begin{equation}
\sin_{\lambda}\big((2m+1)x:a\big)=(2m+1)\sin_{\lambda}\big(x:a\big) \prod_{k=1}^{m}\bigg(1-\cfrac{\sin_{\lambda}^{2}(x:a)}{\sin^{2}\big(\frac{k \pi}{2m+1} \big)}\bigg). \label{28}
\end{equation}
\end{theorem}
Taking $\lambda \rightarrow 0$ in \eqref{28}, we obtain
\begin{equation*}
\sin\big((2m+1)ax\big)=(2m+1)\sin(ax)\prod_{k=1}^{m}\bigg(1-\cfrac{\sin^{2}(ax)}{\sin^{2}\big(\frac{k \pi}{2m+1} \big)}\bigg).
\end{equation*}

\section{Conclusion}
The degenerate trigonometric and degenerate hyperbolic functions are obtained by replacing the ordinary exponentials by the degenerate exponentials in the definitions of trigonometic and hyperbolic functions. In this paper, we considered the degenerate sine, degenerate cosine, degenerate tangent and degenerate cotangent functions and derived several identities among them. Especially, we were able to show the degenerate versions of the following multiple angle formulas (see \eqref{11}, \eqref{15}, \eqref{28}):

\begin{align*}
&2m\cot\big(2mt\big)=-\sum_{j=0}^{2m-1}\tan \Big(t+\frac{j \pi}{2m}\Big), \\
&2m\cot\big(2mt\big)=\sum_{j=0}^{2m-1}\cot(t+\frac{j \pi}{2m}), \\
&\sin\big((2m+1)x\big)=(2m+1)\sin(x)\prod_{k=1}^{m}\bigg(1-\cfrac{\sin^{2}(x)}{\sin^{2}\big(\frac{k \pi}{2m+1} \big)}\bigg).
\end{align*}

\noindent{\textbf{Acknowledgment:}}This paper was supported by Kwangwoon University research grant 2024.

\end{document}